\documentclass[11pt,amstex]{article}

\usepackage{amssymb}
\usepackage{amsmath}
\usepackage{amsfonts}
\usepackage{amscd}
\usepackage{amsthm}

\usepackage[dvips]{color}

\input xypic
\xyoption{all}

\def\a{\alpha}

\def\la{\lambda}

\def\be{\begin{equation}}
\def\ee{\end{equation}}
\def\bear{\begin{eqnarray}}
\def\eear{\end{eqnarray}}
\def\best{\begin{eqnarray*}}
\def\eest{\end{eqnarray*}}

\newtheorem{theorem}{Theorem}[section]

\newtheorem{lemma}[theorem]{Lemma}

\newtheorem{remark}[theorem]{Remark}
\newenvironment{rem}{\begin{remark}\rm}{\end{remark}}
\newtheorem{example}[theorem]{Example}

\newtheorem*{mthm}{Main Theorem}

\def\non{\noindent}
\def\pf{\non {\bf Proof. }}
\def\qed{\nopagebreak \hskip .1in { $\Box$ }\penalty10000 \hskip\parfillskip \par  }
\def\ra{\rightarrow}

\def\del{\overline \partial}

\def\dim{\mbox{\rm dim\,}}

\def\ker{\mbox{ker\,}}

\def\Z{{\mathbb Z}}

\def\cx{{\mathbb C}}

\def\im{\mbox{\rm Im}}
\def\Aut{\mbox{Aut}}

\def\O{{\mathcal O}}

\def\M{{\mathcal M}}

\newcommand{\CM}{\overline{{\mathcal M}}}

\def\im{\mbox{Im}}

\def\la{\langle}
\def\ra{\rangle}

\def\vir{vir}

\def\tf{\tilde{f}}
\def\tJ{\tilde{J}}
\def\tnabla{\tilde{\nabla}}

\def\tomega{\tilde{\omega}}
\def\tg{\tilde{g}}
\def\tL{\tilde{L}}

\addtocounter{section}{0}
\begin{document}

\title{\bf Local  GW Invariants of Elliptic Multiple Fibers}
\vskip.2in

\author{ Junho Lee \\ University of Central Florida \\ Orlando, FL 32816
}

\date{\empty}

\maketitle

\begin{abstract}
We use simple geometric arguments to calculate the dimension zero local Gromov-Witten invariants of elliptic multiple fibers.
This completes the calculation of all dimension zero GW invariants of elliptic surfaces with $p_g>0$.
\end{abstract}

\vspace{1cm}

Let $X$ be a K\"{a}hler surface with $p_g>0$.  By the Enriques-Kodaira classification (cf. \cite{BHPV}),
its minimal model is a $K3$ or Abelian surface, a surface of general type or an elliptic surface.
Each holomorphic 2-form $\a$ on $X$  defines an almost complex structure
\begin{equation}\label{Ja}
J_\a\,=\,(Id+JK_\a)\,J\,(Id+JK_\a).
\end{equation}
Here, $J$ is the complex structure on $X$ and
the endomorphism $K_\a$ of $TX$ is defined by the formula $\la u, K_\a v\ra = \a(u,v)$
where $\la\ ,\ \ra$ is the K\"{a}hler metric. This $J_\a$ satisfies\,:

\begin{lemma}[\cite{L}]\label{ILL}
If $f$ is a $J_\a$-holomorphic map  that represents a nontrivial (1,1) class then its image
lies in the support of the zero divisor $D_\a$ of $\a$ and $f$ is, in fact, holomorphic.
\end{lemma}

The Gromov-Witten invariant $GW_{g,n}(X,A)$ is a (virtual) count of holomorphic maps representing the class $A$.
In particular, the invariant  $GW_{g,n}(X,A)$ vanishes unless $A$ is a (1,1) class since
every holomorphic map represents (1,1) class.
Note that each canonical divisor $D$ of $X$ is a zero divisor of a holomorphic 2-form.
Lemma~\ref{ILL} thus shows that the GW invariant is a sum
$$
GW_{g,n}(X,A)\,=\,\sum GW^{loc}_{g,n}(D_k,A_k)
$$
over the connected components $D_k$ of the canonical divisor $D$ of local invariants
that counts the contribution
of maps whose image lies in $D_k$ (cf. \cite{LP}, \cite{KL}).
It follows that the GW invariants of minimal $K3$ or Abelian surfaces are trivial except possibly for the trivial
homology class because their canonical divisors are trivial.

The local GW invariants have a  universal property.
If $X$ is a minimal surface of general type with a smooth canonical divisor $D$ then
the local invariants associated with $D$, and hence GW invariants, are determined by the normal bundle of $D$ ---
in fact, there exists a universal function of $c_1^2$ and $c_2$ that gives the GW invariants of $X$
(cf. Section 7 of \cite{LP}).

If $\pi:X\to C$ is  a minimal  elliptic surface with $p_g>0$, after suitable deformation, we can assume
$X$ has a canonical divisor of the form
\begin{equation*}\label{can-divisor}
\sum_i n_iF^i\,+\, \sum_k (m_k-1)F_{m_k}
\end{equation*}
where $F^i$ is a regular fiber  and $F_{m_k}$ is a smooth multiple
fiber of multiplicity $m_k$ (cf. Proposition 6.1 of \cite{LP}).
In this case, the GW invariants of $X$ are sums of
universal functions, and are completely determined by the multiplicities $m_k$ and the number
$$c_\pi\,=\,\chi(\O_X)-2\chi(\O_C)$$
(cf. Section 6 of \cite{LP}).  In particular, the
generating function for the set of all {\em dimension zero} GW invariants of $X$  is
given by
\begin{equation}\label{pes}
GW_X^0\,=\,
c_\pi\,
\sum_{d>0}\,  GW_1^{loc}(F,d)\, t^d\, +\,
\sum_k\sum_{d>0}\, GW_1^{loc}(F_{m_k},d) \,t_{m_k}^d
\end{equation}
where the formal variables $t$ and $t_{m_k}$ are
for the fiber class $[F]$ and  the multiple fiber classes $[F_{m_k}]$ respectively; these satisfy
$t_{m_k}^{m_k}=t$. The local invariants in (\ref{pes}) are counts of multiple covers of elliptic curves
together with signs determined by the GW theory of 4-manifolds.

Some of the generating functions  in (\ref {pes}) are known.  In cases of the regular fiber $F$ and the multiple fiber $F_2$, it was proved in Section 10 of \cite{LP} that
\begin{equation}\label{com-known}
GW_1^{loc}(F,d)\,=\,-\frac1d\,\sigma(d)
\ \ \ \   \mbox{and}\ \ \ \
GW_1^{loc}(F_2,d)\,=\,\frac1d\,\left(\sigma(d)\,-\,2\,\sigma\Big(\frac{d}{2}\Big)\right)
\end{equation}
where $\sigma(d)=\sum_{k|d} k$  if $d$ is a positive integer and $\sigma(d)=0$  otherwise.
In this note we use geometric arguments to obtain the terms in (\ref {pes})
associated with fibers of higher multiplicity.
Our main theorem is the following formula for the local invariants $GW_1^{loc}(F_m,d)$ for $m>2$.
This  completes the calculation of all  dimension zero GW invariants of all minimal  elliptic surfaces with $p_g>0$.

\begin{mthm} Let $m\geq 3$. Then
$$
GW_1^{loc}(F_m,d)\,=\,
\frac1d\,\left(\sigma(d)\,-\,m\,\sigma\Big(\frac{d}{m}\Big)\right).
$$
\end{mthm}

\medskip
The contribution of each degree $d$ cover $f$ of elliptic curve $F_m$ is, as a map into a 4-manifold, determined by the normal
bundle $N_m$ of $F_m$.
In cases of $F_1=F$ and $F_2$, the almost complex structure $J_\a$ on $X$ is generic in the sense that
the linearized operator $L_f$ (see (\ref{no}) below) is invertible and hence
the contribution of $f$ is  $(-1)^{h^0(N_m)}/|\Aut(f)|$ (cf. Section 10 of \cite{LP}).
When $m\geq 3$, $J_\a$ is, in general,  no longer generic. We need to perturb $J_\a$ to generic $J$.
In Section 2, using the universal property of local invariants (see (\ref{univ-prop}) below),
we choose a local model that is convenient for our calculation.
In Section 3, when $L_f$ is not invertible, we use a lifting property of covering space to calculate the contribution of $f$
that proves the Main Theorem.
The information for dimension zero
GW invariants  of elliptic surfaces with $p_g>0$ is the same as for its Seiberg-Witten invariants.
We spell out the specific connection in Remark~\ref{SW}.

\bigskip
\noindent
{\bf Acknowledgments.}
I am very grateful to Thomas H. Parker for  helpful discussions.

\vskip 1cm

\setcounter{equation}{0}
\section{Dimension Zero Genus One Local GW Invariants}
\label{section1}
\bigskip

Let  $X$ be a (not necessarily compact) elliptic K\"{a}hler surface with a holomorphic 2-form $\a$.
The 2-form $\a$ defines an almost complex structure $J_\a$ on $X$ by the formula (\ref{Ja}). Suppose that
the zero divisor $D_\a$ of $\a$ has a smooth reduction
$D_\a=(m-1)D$  where $D$ is a regular fiber or a multiple fiber of multiplicity $m$ for some integer $m>1$.
The adjunction formula then shows $c_1(TX)([D])=0$ and $c_1(N)=0$ where $N$ is the normal bundle of $D$.
The moduli space
\begin{equation}\label{l-moduli}
\CM_1^\a(X,d[D])
\end{equation}
of stable $J_\a$-holomorphic maps from  curves of genus one representing
the class $d[D]$ ($d\ne 0$) carries a (virtual) fundamental class
\begin{equation}\label{VFC}
[\CM_1^\a(X,d[D])]^{vir}
\end{equation}
that is defined by the GW theory of 4-manifolds (cf. Section 4 of \cite{LP}).
This (virtual) fundamental class (\ref{VFC}) has dimension zero since $c_1(TX)([D])=0$.
{\em The dimension zero genus one local Gromov-Witten invariant  of $X$ associated with
the zero divisor $D_\a$} is then
\begin{equation*}\label{LGW}
GW_1^{loc}(X,D_\a,d)\,:=\,[\CM_1^\a(X,d[D])]^{vir}.
\end{equation*}
This local GW invariant  has the following universal property.
Let $X^\prime$ be another elliptic K\"{a}hler surface with a holomorphic 2-form $\a^\prime$
whose zero divisor $D_{\a^\prime}=(m^\prime-1) D^\prime$ where $D^\prime$ is a regular fiber or a multiple fiber
of multiplicity $m^\prime$.
Let $N^\prime$ be the normal bundle of  $D^\prime$.
If $m=m^\prime$ and $h^0(N) = h^0(N^\prime)$ then
\begin{equation}\label{univ-prop}
GW_1^{loc}(X,D_\a,d)\,=\,GW_1^{loc}(X^\prime, D_{\a^\prime},d)
\end{equation}
(cf. Section 6 of \cite{LP}). We set
\begin{equation}\label{LGW-g=1}
GW_1^{loc}(X,D_\a,d)\,=\,\left\{
\begin{array}{ll}
GW_1^{loc}((m-1)F,d)\   &\mbox{if \ $D$\ is\ a\ regular\ fiber} \\
GW_1^{loc}(F_m,d)\
&\mbox{if \ $D$\ is\ a\ $m$-multiple\ fiber$^{\displaystyle\phantom{\sum}}$}
\end{array}
\right.
\end{equation}
It was proved in Example 4.4 of \cite{LP} that
\begin{equation}\label{quote}
GW_1^{loc}(mF,d)\,=\,
m\,GW_1^{loc}(F,d).
\end{equation}
As given in (\ref{pes}), all dimension zero GW invariants of minimal elliptic surfaces with $p_g>0$
are sums of local invariants in (\ref{LGW-g=1}).

\bigskip
In the below, we will give a precise description on the (virtual) fundamental class (\ref{VFC})
which will be used for our calculation in Section 3.
The point in the moduli space (\ref{l-moduli})  is an equivalence class  $[f,C]$ of stable maps  $(f,C)$ where
two stable maps $(f,C)$ and $(f^\prime,C^\prime)$
are equivalent if there is a  biholomorphic map $\sigma:C\to C^\prime$ with $f^\prime\circ \sigma=f$.
By Lemma~\ref{ILL}, if $d\ne 0$ every representative $(f,C)$ of  $[f,C]$  is  a holomorphic $d$-fold covering map from $C$ to $D$. Thus,
if $D$ is given by a lattice $\Lambda$ in the complex plane
then $[f,C]$ is determined by an index $d$ sublattice of $\Lambda$.
In particular, the moduli space (\ref{l-moduli}) consists of $\sigma(d)$ points.

On the other hand, since the (virtual) fundamental class (\ref{VFC}) is defined by the GW theory of 4-manifolds,
as described in Section 3 of \cite{IP}, it is a finite sum
\begin{equation}\label{l-sum}
[\CM_1^\a(X,d[D])]^{\vir}\,=\,\sum\,c([f,C])
\end{equation}
over $[f,C]\in\CM_1^\a(X,d[D])$ of the contributions $c([f,C])$ that are defined as follows.
Choose a  $p\in D$ and a small disk $B$ in $X$ with $B\cap D=\{p\}$ and, once and for all,
fix a map $(f,C,x)$ with $f(x)=p$ such that $(f,C)$ represents $[f,C]$.
Then for a generic almost complex structure $J$
on $X$ that is sufficiently close to $J_\a$ and tamed by the K\"{a}hler form on $X$,
there are finitely many $J$-holomorphic maps $(f_i,C_i,x_i)$ from
smooth genus one curves with one marked point such that (i) $f_i(x_i)\in B$
(ii) each  $(f_i,C_i,x_i)$ is $C^0$-close to  $(f,C,x)$ (in a suitable space of maps) and (iii)
the index zero operator
\begin{equation}\label{no}
L_{f_i}\,:\,\Omega^0(f_i^*N_i)\,\to\, \Omega^{0,1}(f_i^*N_i)
\end{equation}
has trivial kernel (or equivalently $L_{f_i}$ is invertible)  where
the operator $L_{f_i}$ is obtained by linearizing $J$-holomorphic map equation (see Remark~\ref{Df} below)
and restricting to the normal bundle $N_i$  of the image of $f_i$.
Denote by
\begin{equation*}\label{one-marked}
\M_{(f,C,x),B,J}
\end{equation*}
the set of such $J$-holomorphic maps $(f_i,C_i,x_i)$.
Notice that for each $(f_i,C_i,x_i)$
the preimage $f_i^{-1}(B)$ consists of $d=|\Aut(f)|$ distinct points $x_{ij}$.
Since the automorphism group of $C_i$ acts transitively, for each $x_{ij}$ there exists an automorphism
$\sigma_j$ of $C_i$ with $\sigma_j(x_i)=x_{ij}$ such that $(f_i\circ \sigma_j,C_i,x_i)$ is also contained in the set
$\M_{(f,C,x),B,J}$. The contribution $c([f,C])$ is thus the (weighted) sum
\begin{equation*}
c([f,C])\,=\,\frac1d\,\sum\, (-1)^{SF(L_{f_i})}
\end{equation*}
over $f_i$ in $\M_{(f,C,x),B,J}$
where the sign of each $f_i$ is given by the mod 2 spectral flow $SF(L_{f_i})$  of the invertible operator $L_{f_i}$.
In particular,  $SF(L_{f_i})=0$ if $L_{f_i}$ is complex linear, namely $J$-linear.

\begin{rem}\label{Df}
The operator $D_{f_i}:\Omega^0(f_i^*TX)\to \Omega^{0,1}(f_i^*TX)$ obtained by
linearizing $J$-holomorphic map equation at $f_i$  is given by
\begin{equation}\label{fDf}
D_{f_i}(\xi)(v)\,=\,\nabla_v\xi +J\nabla_{jv} \xi \,+\,\frac12\big[(\nabla_\xi J) (df_i(jv)) - J(\nabla_\xi J)(v)\big]
\end{equation}
where $\xi\in \Omega^0(f_i^*TX)$, $v\in TC_i$ and $j$ is the complex structure on $C_i$. Here $\nabla$
is the pull-back connection on $f_i^*TX$ of the Levi-Civita connection of the metric on $X$
that is defined by the K\"{a}ler form and $J$ (cf. Lemma 6.3 of \cite{RT}).
\end{rem}

\vskip 1cm

\setcounter{equation}{0}
\section{Local Model}
\label{section2}
\bigskip

Once and for all, fix an integer $m\geq 2$ and let $D$ denote the elliptic curve
given as the complex plane (with coordinate $z$) modulo
the lattice $\Z + (mi)\Z$.
Then $S=D\times \cx$  has an automorphism $\varphi$ of order $m$ defined by
$$
\varphi(z,w)\, =\, \big(z+i, e^{2\pi i/m}\cdot w)
$$
such that all powers $\varphi^i$ are fixed-point free where
$w$ is a coordinate on $\cx$.
Let $S_m$ be the quotient of $S$ by the group $\{\varphi^i\}$
and
$q:S\to S_m$ the quotient map. The map
$S\to\cx : (z,w) \to w^m$
then factors through $S_m$ to give
an elliptic fibration $S_m\to \cx$ whose central fiber is a $m$-multiple fiber $D_m$ given by
the lattice $\Z+i\Z$ with torsion normal bundle $N_m$ of order $m$\,:
\begin{equation*}\label{LocalModel}
\xymatrix{
S=D\times \cx     \ar[rr]^{q}   \ar[dr]_{(z,w)\to w^m}
 &&  \   S_m   \ar[dl] \\
 &   \cx      }
\end{equation*}

\smallskip
The following simple observation is a key fact for our subsequent discussions.
Let $f:C\to D_m$ be a holomorphic map of degree $d$ from an elliptic curve $C$ that is given by
a sublattice of $\Z+i\Z$ of the form
\begin{equation*}\label{lattice}
a\Z \,+\, (bi+k)\Z\ \ \ \ \ \ \ \ \mbox{with}\ \ \ \ \ \ \ d\,=\,ab,\ \  0\,\leq \,k\,\leq \,a-1.
\end{equation*}
Write $D\times \{0\}\subset S$ simply as $D$.

\begin{lemma}\label{key}
Let $D$, $N_m$ and $f:C\to D_m$ be as above. Then,
$$
\mbox{$f$\  factors\  through\  $D$}\ \ \ \Longleftrightarrow\ \ \
a\,|\,\tfrac{d}{m}\ \ \ \Longleftrightarrow\ \ \
f^*N_m\,=\,\O_C
$$
\end{lemma}

\pf
$f$ factors through $D$ $\Longleftrightarrow$ $a\Z+(bi+k)\Z$ is a sublattice of $\Z+(mi)\Z$
$\Longleftrightarrow$ $m|b$ $\Longleftrightarrow$ $a|\frac{d}{m}$. This shows the first assertion.
Observe that for the restriction map $g_m=q_{|D}:D\to D_m$,
\begin{equation}\label{trivial}
g_m^*(N_m)\,=\,g^*_m([D_m]_{|D_m})\,=\,q^*([D_m])_{|D}\,=\,[q^*D_m]_{|D}\,=\,[D]_{|D}\,=\,\O_D
\end{equation}
where $[D_m]$ is the line bundle associated to the divisor $D_m$,
$N_m=[D_m]_{|D_m}$ by adjunction, the pullback divisor $q^*D_m=D$ and
again by adjunction $[D]_{|D}$ is the normal bundle of $D$ that is trivial.
Write as $N_m=\O_{D_m}(p-q)$ where
$$
\int_q^p dz\,=\,\frac{k_1}{m} + i\,\frac{k_2}{m}\ \ \ \ \
\mbox{for\  some}\ \ \ \  0\leq k_1,k_2\leq m-1.
$$
Then, by (\ref{trivial}) and the Abel's Theorem,
$g_m^*N_m=\O_D(\,\sum_j(p_j-q_j)\,)$ for some $p_j,q_j$ such that
$$
\sum_j \int_{q_j}^{p_j} dz\,=\,k_1 + ik_2\,\equiv\,0\ \ \ \mbox{mod}\ \ \ \Z+(mi)\Z.
$$
Consequently, $k_2=0$ and $\gcd(m,k_1)=1$ since $N_m$ is torsion of order $m$.
Now, again by the Abel's Theorem, $f^*N_m=\O_C(\,\sum_\ell(t_\ell-s_\ell)\,)$ for some $s_\ell,t_\ell$ such that
$$
\sum_\ell\int_{s_\ell}^{t_\ell}dz\,=\,
\frac{dk_1}{m}\,\equiv\,0\ \ \ \mbox{mod}\ \ \ a\Z+(bi+k)\Z \ \ \ \Longleftrightarrow\ \ \  f^*N_m=\O_C
$$
Therefore, $a|\frac{d}{m}$ $\Longleftrightarrow$ $m|b$ $\Longleftrightarrow$ $a\,|\,\tfrac{dk_1}{m}$
 $\Longleftrightarrow$ $f^*N_m=\O_C$. This shows the second assertion. \qed

\medskip
\begin{rem}\label{Kform}
Since $q:S\to S_m$ is a covering map, Lemma~\ref{key} shows $f:C\to D_m\subset S_m$ lifts to $\tf:C\to D\subset S$
if and only if $f^*N_m=\O_C$. On the other hand, the K\"{a}hler form $-\frac{i}{2}(dz\wedge d\bar{z}+dw\wedge d\bar{w})$ on $\cx^2$ descends to a
K\"{a}hler form $\tomega$ on $S$ that is $\varphi$-invariant, so $\tomega$ also descends to a K\"{a}hler form
$\omega$ on $S_m$ such that $q^*\omega=\tomega$.
\end{rem}

\vskip 1cm
\setcounter{equation}{0}
\section{Calculation}
\label{section4}
\bigskip

Let $q:(S,D)\to (S_m,D_m)$ be as in Section 2. Fix a holomorphic 2-from
$$
\a\,=\,w^{m-1}dw\wedge dz
$$
on $S$ whose zero divisor is $(m-1)D$ and let
$J_\a$ denote the almost complex structure on $S$ defined by the formula (\ref{Ja}).
The 2-form $\a$ is $\varphi$-invariant, so it descends to a holomorphic 2-form
$\a_m$ on $S_m$ whose zero divisor is $(m-1)D_m$. We denote by $J_m=J_{\a_m}$ the almost complex structure on
$S_m$ defined by the 2-form $\a_m$. Since $D_m$ is a multiple fiber of multiplicity $m$,
\begin{equation*}
GW_1^{loc}(F_m,d)\,=\,[\,\CM_1^{\a_m}(S_m,d[D_m])\,]^{\vir}
\end{equation*}
where the right-hand side is given by the sum of contributions as in (\ref{l-sum}). In order to calculate them,
we decompose the moduli space $\CM_1^{\a_m}(S_m,d[D_m])$ as a disjoint union
$$\ \ \ \ \ \ \
\CM_1^{\a_m}(S_m,d[D_m])\,=\,
\M_{m,d}^+\,{\textstyle \coprod}\,\M_{m,d}^-
\ \ \ \   \mbox{where}\ \ \ \
\left\{
\begin{array}{ll}
[f,C]\in \M_{m,d}^+ &\mbox{if}\ \  h^0(f^*N_m)=0 \\ \,\!
[f,C]\in \M_{m,d}^- &\mbox{if}\ \  h^0(f^*N_m)=1 ^{\phantom{\displaystyle \sum}}
\end{array}
\right.
$$
It then follows from Lemma~\ref{key} that
\begin{equation}\label{cardinality}
\#\,\M_{m,d}^+\,=\,\sigma(d) \,-\,\sigma\Big(\frac{d}{m}\Big)
\ \ \ \ \     \mbox{and}\ \ \ \ \
\#\,\M_{m,d}^-\,=\,\sigma\Big(\frac{d}{m}\Big)
\end{equation}
where $\#A$ is the cardinality of a set $A$.

\medskip

We first calculate the contribution $c([f,C])$ of $[f,C]$ in $\M_{m,d}^+$.
In the below, we always assume $m\geq 3$ and $m|d$.

\begin{lemma}\label{step1}
If $[f,C]\in\M_{m,d}^+$ then $c([f,C])=\frac1d$.
\end{lemma}

\pf
The linearized operator $L_f$ has the form
$L_f=\del_f+R_m$ where $\del_f$ is the usual $\del$-operator on $f^*N_m$ and the zeroth order term $R_m$ is given by
\begin{equation*}\label{bundlemap}
R_m(\xi)\,=\,-\nabla_\xi K_{\a_m}\circ J_{\a_m}\circ df
\ \ \ \ \ \mbox{for}\ \ \ \ \xi\in \Omega^0(f^*N_m)
\end{equation*}
(cf. Section 8 of \cite{LP}). But, $R_m\equiv 0$ since
$\a_m$ (and hence $K_{\a_m}$) vanishes of order $m-1\geq 2$ along $D_m$. Consequently,
$\dim\ker L_f=2h^0(f^*N_m)=0$, so $L_f$ is invertible with $SF(L_f)=0$.
Now, the proof follows from the fact $f:C\to D_m$ has degree $d$.
\qed

\bigskip
Let $[f,C]\in \M_{m,d}^-$. The proof of Lemma~\ref{step1} shows
$L_f=\del_f$ is not invertible. In this case, we will uses the $m$-fold covering map $q:S\to S_m$
to calculate the contribution $c([f,C])$.
Observe that  by Lemma~\ref{key} the map
$$
\M_{m,d}^- \,\to\, \CM^\a_1(S,\tfrac{d}{m}[D])
\ \ \ \ \ \mbox{defined\ by}\ \ \ \ \
[f,C]\to [\tf,C]
$$ is one-to-one and onto where $\tf$ is a lift of $f$.

\begin{lemma}\label{step2}
If $[f,C]\in \M_{m,d}^-$ then $c([f,C])=\tfrac1m\,c([\tf,C])$.
\end{lemma}

\pf
Let $B=\{0\}\times \Delta\subset S$ where $\Delta$ is a small disk around $0$ in $\cx$ and $B_m=q(B)$ and
fix a map $(f,C,x)$ with $f(x)\in B_m$ such that $(f,C)$ represents $[f,C]$.
Since the restriction map $q_{|B}:B\to B_m$ is one-to-one,
Lemma~\ref{key} shows that $(f,C,x)$
uniquely lifts to a $J_\a$-holomorphic map $(\tf,C,x)$ with $\tf(x)\in B$ such that
$(\tf,C)$ represents $[\tf,C]$ in $\CM_1^\a(S,\frac{d}{m}[D])$.

Let $\omega$ and $\tomega$ be the K\"{a}hler forms as in Remark~\ref{Kform} and
choose a generic $\omega$-tamed almost complex structure $J$ on $S_m$ that is close to $J_m$. Then, we have
\begin{itemize}
\item
$J$ lifts to an $\tomega$-tamed almost complex structure $\tJ$ on $S$ close to $J_\a$ such that
$dq\circ \tJ=J\circ dq$,
\item
each $f_i$ in $\M_{(f,C,x),B_m,J}$ is homotopic to $f$ since $f_i$ is $C^0$-close to $f$, so
$(f_i,C_i,x_i)$ also uniquely lifts to $\tJ$-holomorphic maps $(\tf_i,C_i,x_i)$ with $\tf(x_i)\in B$ such that
$(\tf_i,C_i,x_i)$ is $C^0$-close to $(\tf,C,x)$.
\end{itemize}

The pair $(\omega,J)$ defines a metric $g$ on $S_m$
whose lift $\tg=q^*g$  is the same  metric defined by the pair $(\tomega,\tJ)$.
Let $\nabla$ and $\tnabla$ respectively denote the pull-back connections on $f_i^*TS_m$ and  $\tf_i^*TS$
of the Levi-Civita connection of $g$ and  $\tg$.
The differential $dq$ then induces a bundle isomorphism $dq:\tf_i^*TS\to \tf_i^*q^*TS_m=f_i^*TS_m$ such that
$dq \circ \tnabla = \nabla \circ dq$  (see \cite{W}  page 138) and hence by the formula (\ref{fDf}) we have
\begin{equation}\label{variation}
dq\circ D_{\tf_i}\,=\,D_{f_i}\circ dq
\end{equation}
The differential $dq$ also induces a bundle isomorphism $dq_i:\tf_i^*\tilde{N}_i\to f_i^*N_i$ and
restricting the equation (\ref{variation}) to  $\tf_i^*\tilde{N}_i$ and $f_i^*N_i$ gives
$$
dq_i\circ L_{\tf_i}\,=\,L_{f_i}\circ dq_i
$$
where $\tilde{N}_i$ and $N_i$ are normal bundles of $\im(\tf_i)$ and   $\im(f_i)$ respectively.
Therefore, $L_{\tf_{i}}$ is also invertible and hence there is
one-to-one correspondence
$$
\M_{(f,C,x),B_m,J}\,\to \,\M_{(\tf,C,x),B,\tJ}\ \ \ \ \
\mbox{given\ by}\ \ \ \ \
(f_i,C_i,x_i)\,\to\,(\tf_i,C_i,x_i).
$$

Let $\tL_t$ be a path
of first  order elliptic operators from an invertible $\tJ$-linear
operator $\tL_0$ to $\tL_1=L_{\tf_i}$ with all $\tL_t$ invertible except
at finitely many $t_k$. Then, $dq_i\circ \tL_{t}\circ (dq_i)^{-1}$ is also a path from
invertible $J$-linear operator to $L_{f_i}$ such that
$$
SF(L_{\tf_i})\,\equiv\,\sum_k \dim \ker \tL_{t_k}\,=\,
\sum_k \dim \ker dq_i\circ \tL_{t_k}\circ (dq_i)^{-1}\,\equiv\,SF(L_{f_i}) \qquad \mbox{(mod 2)}.
$$
Now, noting $\deg(f)=d$ and $\deg(\tf)=\frac{d}{m}$, we have
\begin{equation*}\label{E3}
c([f,C])\,=\,
\frac1d\,\sum_{f_i}\, (-1)^{SF(L_{f_i})}
\,=\,\frac{1}{d}\,
\sum_{\tf_i}\, (-1)^{SF(L_{\tf_i})}\,\,=\,\frac1m\,c([\tf,C]).
\mbox{\qed}
\end{equation*}

\bigskip
We are now ready to prove the Main Theorem in the introduction.

\bigskip
\non
{\bf Proof of the Main Theorem :}
It follows from Lemma~\ref{step1}, Lemma~\ref{step2} and (\ref{cardinality}) that
\begin{align}\label{last1}
GW_1^{loc}(F_m,d)\,&=\,
\sum_{[f,C]\in\M_{m,d}^+} c([f,C]) \,+\,
\sum_{[f,C]\in\M_{m,d}^-} c([f,C]) \notag \\
&=\,
\frac1d\,\left(\sigma(d) - \sigma\Big(\frac{d}{m}\Big)\right)\,+\,
\frac1m\,[\,\CM^\a_1(S,\tfrac{d}{m}[D])\,]^{\vir}.
\end{align}
Since the 2-form $\a$ on $S$ has the zero divisor $(m-1)D$, so by (\ref{quote}) and (\ref{com-known})
we have
\begin{equation}\label{last2}
[\,\CM^\a_1(S,\tfrac{d}{m}[D])\,]^{\vir}\,=\,GW_1^{loc}((m-1)F,\tfrac{d}{m})\,=\,
-(m-1)\,\frac{m}{d}\,\sigma\Big(\frac{d}{m}\Big).
\end{equation}
Now, the proof follows from (\ref{last1}) and (\ref{last2}). \qed

\bigskip
\begin{rem}
One can also use the above argument  to compute $GW_1^{loc}(F_2,d)$, replacing the ``Taubes type'' argument used in \cite{LP}.
Specifically, for each $f\in\CM_1^{\a_2}(S_2,d[D_2])$ the linearized operator $L_f$ is invertible
with $SF(L_f)\equiv h^0(f^*N_2)$ (mod 2) (cf. Proposition 9.2 of \cite{LP}).
Thus, by (\ref{cardinality}) we have
$$
GW^{loc}_1(F_2,d)\,=\,[\,\CM_1^{\a_2}(S_2,d[D_2])\,]^{\vir}\,=\,
\frac1d\,\left(\sigma(d)-\sigma\Big(\frac{d}{2}\Big)\right)\,-\,\frac1d\,\sigma\Big(\frac{d}{2}\Big)\,=\,
\frac1d\,\left(\sigma(d)-2\sigma\Big(\frac{d}{2}\Big)\right).
$$

\end{rem}

\medskip

\begin{rem}\label{SW}
Ionel and Parker \cite{IP} showed how  GW invariants for the class $A$
of  a symplectic 4-manifold $X$ are related with the
Taubes' Gromov invariants $Gr_X(A)$ \cite{T} that count
embedded (not necessarily connected) curves in  $X$ representing the class $A$.
They used a particular function $F(t)$ that satisfies
$$
\underset{d}{\prod} \,{F\big(t^d\big)}^{ -\frac1d\sigma(d)}\,=\,(1-t)
$$
to relate Taubes' counting of multiple covers of embedded tori with the dimension zero genus one GW invariants.
Let $X$ be a minimal elliptic surface with $p_g>0$. In this case, any GW invariant constrained to pass through
generic points vanishes (cf. Corollary 3.4 of \cite{LP}).
So, by (\ref{pes}), (\ref{com-known}) and the Main Theorem,
the relation between two set of invariants (Theorem 4.5 of \cite{IP}) yields
\begin{align*}
\sum_A\,Gr_X(A)\,t_A\,&=\,
\prod_{d,k}\,F(t^d)^{c_\pi GW_1^{loc}(F,d)}\,F(t_{m_k}^d)^{GW_1^{loc}(F_{m_k},d)}\\
&=\,(1-t)^{c_\pi}\,\prod_{k}\,(1+t_{m_k}+\cdots + t_{m_k}^{m_k-1}).
\end{align*}
This also gives the well-known Seiberg-Witten invariants $SW$ of $X$ (cf. \cite{FM}, \cite{B}, \cite{FS}) due to
the famous Taubes' theorem $SW=Gr$.
\end{rem}


\begin{thebibliography}{}








\bibitem[B]{B} R. Brussee,
{\em The canonical class and the $C^\infty$ properties of Kähler surfaces}, New York J. Math. {\bf 2} (1996), 103--146.


\bibitem[BHPV]{BHPV} W. Barth, K. Hulek, C. Peters, and A. Van de Ven,
{\em Compact complex surfaces},  second  ed., Springer-Verlag, Berlin
Heidelberg, 2004.







\bibitem[FM]{FM} R. Friedman and J. Morgan,
{\em Obstruction bundles, semiregularity, and Seiberg-Witten invariants},
 Comm. Anal. Geom. {\bf 7} (1999), no. 3, 451--495.







\bibitem[FS]{FS} R. Fintushel and R. Stern,
{\em Rational blowdowns of smooth 4-manifolds}, J. Diff. Geom. {\bf 46} (1997), 181-235.




\bibitem[IP]{IP} E. Ionel and T.H.  Parker,
{\em The Gromov invariants of Ruan-Tian and Taubes}, Math. Res. Lett. {\bf 4} (1997), no. 4, 521--532.




\bibitem[KL]{KL} Y-H. Kiem and J. Li,
{\em Low degree GW invariants of spin surfaces}, Pure Appl. Math. Q. {\bf 7} (2011), no. 4, 1449--1476.





\bibitem[L]{L}  J. Lee, \textit{Family Gromov-Witten Invariants for
K\"{a}hler Surfaces},  Duke Math. J. \textbf{ 123} (2004), no. 1, 209--233.





\bibitem[LP]{LP} J. Lee and T.H. Parker,
{\em A Structure Theorem for the Gromov-Witten Invariants of K\"{a}hler Surfaces},
 J. Diff. Geom. {\bf 77} (2007), 483-513.




















\bibitem[RT]{RT} Y. Ruan and G. Tian, {\em A mathematical theory of 
cohomology}, J. Differential Geom. {\bf 42} (1995), 259-367.





\bibitem[T]{T} C. Taubes,
{\em Seiberg Witten and Gromov invariants for symplectic $4$-manifolds},
Edited by Richard Wentworth. First International Press Lecture Series, 2. International Press, Somerville, MA, 2000.


\bibitem[W]{W}  G. Walschap, {\em Metric structures in differential geometry}, Graduate Texts in Mathematics, 224.
Springer-Verlag, New York, 2004.


\end{thebibliography}
\end{document}